\documentclass[twoside]{article}
\usepackage{amssymb,amsmath,epsfig,mathrsfs}                  
\begin{document}                                              
\catcode`!=11                                                 
\let\!int\int \def\int{\displaystyle\!int}                    
\catcode`!=12                                                 
\def\pd#1#2{\frac{\partial#1}{\partial#2}}                    
\baselineskip=5truemm \font\abstractsizebf=cmbx9              
\font\abstractsizeit=cmti9 \font\shumeiit=cmti8               
\font\abstractsize=cmr9 \font\shumeisize=cmr8                 
\abovedisplayskip=.5mm plus .5mm minus .5mm                   
\belowdisplayskip=.5mm plus .5mm minus .5mm                   
\def\evenhead{\pushziti\hbox to                               
\hsize{\rm\shumeisize\thepage\hfill                           
    \shumeisize  \hspace{10mm}COMM. MATH. RES.                
    \hfill VOL.~3?}\popziti}                                  
\def\oddhead{\hbox to \hsize{\shumeisize NO.~?\hfill          
    {\shumeiit XXXX X. X. et al.}\ \
    \shumeisize A THEOREM OF NEHARI TYPE
    \hfill\rm\shumeisize\thepage}}                            
\title{\vspace{-1in}                                          
   \parbox{\linewidth}{\footnotesize\noindent                 
     \unitlength=1mm                                          
     \doublerulesep0.5pt                                      
}                                                             
   \vspace{10mm}
\\                                                            
\bf On the Lower Bound of the Divisibility of Exponential Sums   in Binomial Case                                        
\footnotetext{                                                
\bf\scriptsize Received date: \rm\scriptsize Jan. 26, 2012.   
\hfil\break\indent                                            
\bf\scriptsize\hspace{1.3mm} Foundation item: \rm\scriptsize  
The NSF 61502230 of China.                                  
\hfil\break\indent                                            
\bf\scriptsize\hspace{2.3mm}{\rm *} {Corresponding author{\rm .}} 
\hfil\break\indent                                            
\bf\scriptsize\hspace{2.3mm}{E-mail address: \rm\scriptsize   
liuxg0201@163.com (Liu X)
} 
}}
\date{}
\author{\abstractsize\sc 
Liu Xiaogang$^{  {\text *}}$
\\
\abstractsize{({\it\abstractsizeit   School of Computer Science and Technology , Nanjing Tech University, Nanjing, 211800})}\\  
}

\maketitle


\let\oldparindent\parindent
\leftskip10truemm
\rightskip10truemm
{
\parindent=0pt
\abstractsizebf Abstract: \rm
\abstractsize
In this article, we analyze the lower bound of the divisibility of families of exponential sums for  binomials   over prime field.
An upper bound is given for the lower bound, and, it is related to permutation polynomials.                                                    
\\
\abstractsizebf Key words:
\abstractsize exponential sums, finite fields, binomials, permutation polynomials                  
\\
\abstractsizebf 2010 MR subject classification:
\abstractsize 11L07, 12E30                                   
\\
\abstractsizebf Document code:
\abstractsize A
\\
\abstractsizebf Article ID:
\abstractsize 1674-5647(201?)0?-0???-0?
\\
\abstractsizebf DOI:
\abstractsize 10.13447/j.1674-5647.????.??.??
}

\def\parindent{\oldparindent}
\leftskip0truemm
\rightskip0truemm


\newtheorem{theorem}{Theorem}[section]                                               %
\newtheorem{definition}{Definition}[section]                                         %
\newtheorem{lemma}{Lemma}[section]                                                   %
\newtheorem{example}{Example}[section]
\newtheorem{proposition}{Proposition}[section]                                       %
\newtheorem{corollary}{Corollary}[section]                                           %
\newtheorem{remark}{Remark}[section]                                                 %
\renewcommand{\theequation}{\thesection.\arabic{equation}}                           %
\catcode`@=11                                                                        %
\@addtoreset{equation}{section}                                                      %
\setcounter{page}{1}                                                                 %
\thispagestyle{empty}                                                                %
\catcode`@=12                                                                        %

\section{Introduction}\label{secI}

Exponential sums and their divisibility have been applied to characterize important properties of objects in applied mathematics.  There are many estimates for the divisibility of exponential sums \cite{AS,CFM,MM,MSC,S}.  It is difficult  in general. 
 The exact divisibility of families of exponential sums associated to
binomials $F(X)=aX^{d_1}+bX^{d_2}$, is computed under some natural conditions when $a,b \in \mathbb{F}_p^*$ \cite{CFM}.
The result is applied to solutions of equations,  to the study of Waring problem over finite fields, and to the determination of permutation polynomials.

 Let $F(X)$ be a two term polynomial in one variable over prime field $\mathbb{F}_p$.
The following bound  for the valuation (divisibility) of an exponential sum
can be thought of as a particular case of  \cite[Theorem 8]{MSC}.
\begin{theorem} \label{thm1}
Let  $ F(X)=aX^{d_1}+bX^{d_2}, \mbox{with}\  a,b\in \mathbb{F}_p^*, 1\leq d_1\not= d_2\leq p-2$. If $S_p(F)$  is the exponential sum $ \sum\limits_{X\in \mathbb{F}_p}\phi (F(X))$, then $\nu_\theta (S_p(F))\geq \mu_p (d_1,d_2)$, where 
\[
\mu_p(d_1,d_2)=\min\limits_{(i,j)}\left\{i+j|0\leq i,j <p\right\},
\]
 for $(i, j )\not=(0, 0)$ a solution to the modular equation
\[
d_1i+d_2j \equiv 0 \ \mbox{mod}  \ p-1.
\]
 
\end{theorem}

\begin{remark} \quad %
{\rm
Let $\mathbb{Q}_p$ be the $p$-adic field and let $\xi$ be a primitive $p$-th root of unity in $\overline{\mathbb{Q}}_p$. Define $\theta =1-\xi$ and denote by $\nu_{\theta}$ the valuation over $\theta$. Note that $\nu_{\theta} (p)=p-1$ and $\nu_p (x)={{\nu_{\theta}(x)}\over {p-1}}$. Let $\phi :\mathbb{F}_p \to \mathbb{Q}(\xi)$ be the nontrivial additive character defined by $\phi (a)=\xi^{a}$, for $a\in \mathbb{F}_p$. The exponential sum associated to $F(X) $ is defined as 
$S_p(F)=\sum\limits_{X\in \mathbb{F}_p}\phi (F(X)). $
}
\end{remark}

In the next section, Theorem \ref{MR001} is presented on the lower bound  in Theorem \ref{thm1}, which is less than half of $p$. To this end,   the domain $[{0}, 1)$ is splitted into infinitely many segments,  which are for different range of integers. In this process,  we can make a success in our analysis for every step, which is not enough for the whole.  But finally, with the splitting parts changing and interacting with each other, we can make it complete. Though not much long and complicated, the important idea in our work is to note the effect of those seemingly nonsignificant interactions between different segments. 

\section{Main Result}\label{SecIII}

 \begin{theorem}\label{MR001}
Let $1\leq d_1\not= d_2\leq p-2$ be positive integers and let $p\geq 5$ be a prime. Then \begin{equation}\label{R001}
  \mu_p(d_1,d_2)\leq {{p-1}\over 2}.
\end{equation}
.
\end{theorem}
\noindent{\it Proof.}
First, consider the case when $d_1,d_2$ are odd numbers, and $\mbox{gcd}(d_1,p-1)=\mbox{gcd}(d_2,p-1)=1$. 

Let $i=1$. Then $d_1i+d_2j=d_1+d_2j$.  If it is zero modular $p-1$, then $j\not=p-1$; if $j=p-2$ then $d_1+d_2j=d_1+d_2(p-2)=d_1+d_2(p-1)-d_2$ is not equal to zero modular $p-1$. Also, $j\not=0$.  Note that such a $j$ must exist, since $\mbox{gcd}(d_2,p-1)=1$. So,  
\begin{equation}\label{R004}
1\leq j\leq p-3.
\end{equation} 

If $j>{{p-1}\over 2}$,   consider the following sets
\[
\left[
1-{1\over {2^{k-1}}},1-{1\over {2^{k }}}
\right)
\]
for $k=1,2,3,\ldots.$ We can find that the union of the sets is $[{0}, 1)$. Assume that
$j=\alpha_{k-1} (p-1)$ and $\alpha_{k-1} \in\left[
1-{1\over {2^{k-1}}},1-{1\over {2^{k }}}
\right) $. Then 
\[
1-{1\over {2^{k-1 }}} \leq \alpha_{k-1}
\]
which implies that
\[
\left(
1-{1\over {2^{k-1 }}}\right) (p-1)\leq p-3 (\mbox{equ:} \ref{R004})
\]
that is
\[
p-1-(p-3) \leq{{p-1}\over {2^{k-1}}} . 
\]
Thus
\begin{equation}\label{R002}
{{p-1}\over {2}}\geq 2^{k-1}.  
\end{equation}

Since $d_1+d_2j \equiv 0 \ \mbox{mod}  \ p-1$, we have
\[
d_1\cdot 2+d_2\cdot 2j \equiv 0 \ \mbox{mod}  \ p-1,
\]
and, 
\[
\begin{array}{lll}
2j\in  2\left[
1-{1\over {2^{k-1}}},1-{1\over {2^{k }}}
\right)(p-1) 
&=&\left[
2-{1\over {2^{k-2}}},2-{1\over {2^{k-1}}}
\right)(p-1)\\
&=&p-1+\left[
1-{1\over {2^{k-2}}},1-{1\over {2^{k-1}}}
\right)(p-1).
\end{array}
\]
So,  
\[
2j \equiv j_{k-2} \ \mbox{mod}  \ p-1
\]
for $ j_{k-2}
   \in \left[
1-{1\over {2^{k-2}}},1-{1\over {2^{k-1}}}
\right)(p-1)$, and
\[
d_1\cdot 2+d_2\cdot j_{k-2} \equiv 0 \ \mbox{mod}  \ p-1.
\]
Multiply both sides of the modular equation by 2 again, and continue in this way, we find that
\[
d_1\cdot 2^{k-1}+d_2\cdot j_{0} \equiv 0 \ \mbox{mod}  \ p-1
\]
for $ j_{0}
   \in \left[
 0,{1\over 2}
\right)(p-1)$. Denote $i'= 2^{k-1} \in\left(
 0,{{p-1}\over 2}
\right]$ (equ: \ref{R002}), and $j'=j_0$. If $i'+j'\leq {{p-1}\over 2}$, we are done.
Otherwise,  let 
\[
i''={{p-1}\over 2}-i',  \  \ j''={{p-1}\over 2}-j'.
\]
Then
\[
\begin{array}{lll}
&&d_1\cdot i''+d_2\cdot j''= d_1\cdot \left({{p-1}\over 2}-i'\right)+d_2\left({{p-1}\over 2}-j'\right)=(d_1+d_2)\left({{p-1}\over 2}\right)-(d_1\cdot i'+d_2\cdot j')\\
& \equiv & -(d_1\cdot i'+d_2\cdot j') \ \mbox{mod}  \ p-1\\
& \equiv & 0 \ \mbox{mod}  \ p-1
\end{array}
\]
 we have that $d_1,\ d_2$ are odd numbers and $d_1+d_2$ is even. And 
\[
i''+j''={{p-1}\over 2}-i'+{{p-1}\over 2}-j'=p-1-(i'+j')\leq {{p-1}\over 2}.
\] 

Now,  assume that $g_1=\mbox{gcd}(d_1,p-1)  \geq 2$ .
Let $i={{p-1}\over g_1}$ and $ j=0$. Then
 \[
d_1\cdot i +d_2\cdot j \equiv 0 \ \mbox{mod} \ p-1,
\] 
and
 \[
i+j={{p-1}\over g_1}\leq {{p-1}\over 2}.  
\]
This completes the proof of the theorem.    \quad

\begin{remark} \quad %
{\rm
Theorem \ref{MR001} is closely related to permutation polynomials.  A permutation polynomial over a finite field is one which permutes the elements of the field. There is the  Hermite's criterion about whether a polynomial is a permutation polynomial
\begin{theorem}$^{[6,7]}$\label{tm2.2}
A polynomial $F(X)$ over a finite field $\mathbb{F}_q $ ($q$ is a power of prime $p$)  is a permutation polynomial if and only if 

\renewcommand{\labelenumi}{$($\mbox{\arabic{enumi}}$)$}
\begin{enumerate} 
 \item for each $i$ with $0<i<q-1$, the reduction of $F(X)^i$ modulo $X^q-X$ has degree less than $q-1$; and 
 \item $F(X)$ has exactly one root in $\mathbb{F}_q$.
\end{enumerate}

\end{theorem}

In \cite{CFG}, Roger improved this result by showing that only the first ${p-1}\over 2$ reductions are needed to be tested for $F(X)$ to be a permutation polynomial when $\mathbb{F}_q$ is the prime field $\mathbb{F}_p$. 
If there exist $d_1,d_2$, such that $\mu_p(d_1,d_2)>{{p-1}\over 2}$, then Roger's result implies that the     first condition in Theorem \ref{tm2.2} is satisfied when testing whether the binomial $F(X)=aX^{d_1}+bX^{d_2}$ is a permutation polynomial. In fact, for $d \leq {{p-1}\over 2}$, when 
\[
F(X)^d= (aX^{d_1}+bX^{d_2} )^d=\sum\limits_{i=0}^{d}
\left(
\begin{array}{c}
d \\
i \\
\end{array}
\right)
a^ib^{d-i}X^{d_1i+d_2(d-i)}.
\]
is expanded and reduced, there will be no terms of degree $p-1$. To see this, assume
\[
d_1i+d_2(d-i)=c(p-1)+d_0
\]
for some nonnegative integers $c$ and $0 \leq d_0 <p-1$. By assumption, $\mu_p(d_1,d_2)>{{p-1}\over 2}$, now $d=i+(d-i)  \leq {{p-1}\over 2}$, so   $d_0 \not= 0$. Since $X^{p-1}=1$ for $X\not=0$, $F(X)^d$ has no degree $p-1$ terms. 

But our result (equ: \ref{R001}) implies that such a situation will never happen. That is, generally speaking,  the exponents $d_1,d_2$ can't be parameterized to satisfy the first condition of Theorem \ref{tm2.2} using Roger's result. 
}
\end{remark}



Thus, when $\mathbb{F}_q$ is the prime field $\mathbb{F}_p$, there are $p-2$  reductions to be tested for checking whether a polynomial $F(X)$ is a permutation polynomial. In \cite{CFG}, Roger improved this result by showing that only the first ${p-1}\over 2$ reductions are needed to be tested.


\end{document}